\documentclass[12pt,a4paper,reqno]{article}
\def\hybrid{\topmargin 0pt      \oddsidemargin 0pt
        \headheight 0pt \headsep 0pt
        \textwidth 160true mm       
        \textheight 231true mm         
        \marginparwidth 0.0in
        \parskip 3pt plus 1pt   \jot = 1.5ex}

\hybrid
\usepackage{amssymb}
\usepackage{amsmath}
\usepackage{amsthm}

\newcommand{\g}{\mathfrak{g}}

\newcommand{\lf}{\mathfrak l}

\newcommand{\ff}{\varphi}

\newcommand{\Uh}{{U}_h (\g)}
\newcommand{\Ug}{U(\g)}

\newcommand{\A}{\mathcal{A}}

\newcommand{\C}{\mathbb{C}}

\renewcommand{\O}{\mathcal{O}}
\newcommand{\R}{\mathcal{R}}

\newcommand{\ot}{\otimes}

\renewcommand{\[}{[\![}
\renewcommand{\]}{]\!]}

\newcommand{\hlf}[1]{\frac{#1}{2}}

\newcommand{\Vect}{\operatorname{Vect}}
 
\newcommand{\eps}{\varepsilon}

\theoremstyle{plain} 
\newtheorem{thm}{Theorem}[section]
\newtheorem{cor}{Corollary}[section]
\newtheorem{lemma}{Lemma}[section]
\newtheorem{propn}{Proposition}[section]

\theoremstyle{definition}
\newtheorem{defn}{Definition}[section]

\newtheorem{re}{}[section]
\newtheorem{conjecture}[re]{Conjecture}

\theoremstyle{definition}  
\newtheorem{rem}{Remark}[section]

\newcommand{\be}[1]{\begin{eqnarray#1}}
\newcommand{\ee}[1]{\end{eqnarray#1}}

\newcommand{\De}{\Delta}
\newcommand{\de}{\delta}

\renewcommand{\t}{\otimes}

\newcommand{\bref}[1]{(\ref{#1})}
\newcommand{\tr}{\triangleright}         

\newcommand{\Sym}{\operatorname{Sym}}
\numberwithin{equation}{section}

\title{Quantum $G$-manifolds}
\author{J. Donin\thanks{This research is partially supported
by Israel Academy of Sciences Grant no. 8007/99-01}\\
{\normalsize Dept. of Math.  Bar-Ilan University}\\
{\normalsize Max-Planck-Institut f\"ur Mathematik}}

\begin{document}
\date{}
\maketitle
\begin{abstract}
Let $G$ be a Lie group, $\g$ its Lie algebra, and $U_h(\g)$ the corresponding
quantum group. 
We consider some examples of $U_h(\g)$-invariant one and two parameter quantizations
on $G$-manifolds.
\end{abstract}
\section{Introduction}

Passing from  classical mechanics to quantum mechanics involves
replacing the commutative function algebra, $\A$, of classical observables
on the appropriate phase space, $M$, with a noncommutative (deformed) algebra,
$\A_t$, of quantum observables (see \cite{BFFLS} where the 
deformation quantization scheme is developed). 
The algebra $\A$ is a Poisson algebra and the product in $\A_t$ 
is given by a power series in the formal parameter $t$ with leading term 
the original commutative product and with leading term in the
commutator   given by the Poisson  bracket.
If the classical system is invariant under a 
Lie group of symmetries, $G$, then $G$ acts on the phase space (and thus on $\A$) and
the associated quantum system often retains 
the group of symmetries. In particular, the algebra $\A_t$ is often 
invariant under the action of $G$, or under the action of its universal 
enveloping algebra $U(\g)$.

Modern field-theoretical models and, in particular, the problem of 
incorporating gravity into a quantum field theory led to the requirement
of deforming  (quantizing) the group symmetry and the phase space themselves.
This is one of the reasons for the interest in quantum groups. 
The quantum group, $\Uh$, defined by Drinfeld and Jimbo is a deformation of 
$U(\g)$ as a Hopf algebra.
The quantization of the phase space and its  symmetry group corresponds to
a $\Uh$ invariant deformation of the  algebra $\A_t$, which leads us
to the problem of two parameter (or double) quantization, $\A_{t,h}$, 
of the function algebra for a $U(\g)$ invariant Poisson structure.
In other words, the problem of two parameter (or double) quantization appears if
we want to quantize the Poisson bracket in such a way that 
multiplication in the quantized algebra is invariant under the 
quantum group action.

In the present talk we consider examples of one and two parameter
invariant quantizations of Poisson function algebras on some natural $G$-manifolds. 

In fact, we will regard  $\A$ as the sheaf of algebras and $\A_t$ ( $\A_{t,h}$ )
as a deformation of that sheaf. We call the {\em quantum} $G$-manifold
the manifold $M$ endowed with the deformed sheaf of non-commutative algebras.
  
In Section 2 we review some facts on quantum groups related to our problem.
Note that we consider only the case of semisimple $G$.
In Section 3 we formulate precisely the problem and reduce it to
the problem of a $G$-invariant non-associative quantization. Note that 
in our setting we suppose that the action of $G$ on the space of function on $M$
does not deform. In Section 4 we describe Poisson brackets admissible
for $\Uh$-invariant one and two parameter quantizations.

Beginning from Section 5 we consider examples of $\Uh$-invariant quantization
on some $G$-manifolds.
Namely, we consider two-sided and $Ad$-action $G$ on itself, coadjoint
action $G$ on $\g^*$, and  semisimple orbits in $\g^*$. 

{\bf Acknowledgment.} I thank Max-Planck-Institut f\"ur Mathematik for
hospitality and very stimulating working atmosphere.

\section{Quantum groups}

We will consider quantum groups in the sense of Drinfeld, \cite{Dr2},
as Hopf algebras being deformed universal enveloping algebras. So we will regard that
$\Uh=U(\g)[[h]]$ as a topological $\C[[h]]$-module and
$\Uh/h\Uh=U(g)$ as a Hopf algebra over $\C$. 
In particular, the deformed comultiplication in  $\Uh$ has the form
\be{}
\label{comult}
\De_h=\De+h\De_1+o(h),
\ee{}
where $\De$ is the comultiplication in the universal enveloping algebra $U(\g)$.
One can prove, \cite{Dr2},  that the map
$\De_1:U(\g)\to U(\g)\ot U(\g)$ is such that $\De_1-\sigma\De_1=\de$
($\sigma$ is the usual permutation) being restricted to $\g$ gives
a map $\de:\g\to\wedge^2\g$ which is a 1-cocycle and defines
the structure of a Lie coalgebra on $\g$ (the structure of a Lie algebra on
the dual space $\g^*$).  
The pair $(\g,\de)$ is called
a quasiclassical limit of $\Uh$.

In general, a pair $(\g,\de)$, where
$\g$ is a Lie algebra and $\de$ is such a 1-cocycle, is called a Lie bialgebra.
It is proven, \cite{EK}, that any Lie bialgebra $(\g,\de)$ can be quantized, i.e., 
there exists a quantum group  $\Uh$ such
that  the  pair  $(\g,\de)$ is its quasiclassical limit.

A Lie bialgebra $(\g,\de)$ is said to be  a coboundary one if  there exists an element
$r\in\wedge^2\g$, called the classical r-matrix, 
such that $\de(x)=[r,\De(x)]$ for $x\in\g$. Since $\de$ defines a Lie coalgebra
structure, 
$r$ has to satisfy the so-called classical Yang-Baxter equation
which can be written in the form
\be{}
\label{cYB}
\[r,r\]=\ff,
\ee{}
where $\[\cdot,\cdot\]$ stands for the Schouten bracket and 
$\ff\in\wedge^3\g$ is an invariant element. We denote the coboundary
Lie bialgebra by $(\g,r)$.

In case $\g$ is a simple Lie algebra, the most known r-matrix is
the Sklyanin-Drinfeld one:
\be{}\label{SDr}
r=\sum_\alpha X_\alpha\wedge X_{-\alpha},
\ee{} 
where the sum runs over all positive roots; the root vectors $X_\alpha$
are chosen in such a way that $(X_\alpha, X_{-\alpha})=1$
for the Killing form $(\cdot,\cdot)$.
This is the only r-matrix of weight zero, \cite{SS}, and its quantization is the Drinfeld-Jimbo
quantum group.
A classification of all r-matrices for simple Lie algebras was
given in \cite{BD}.

From results of Drinfeld and of Etingof and Kazhdan one can derive
the following description of quantum groups associated with
semisimple Lie algebras.
\begin{propn}\label{propo2.1}
Let $\g$ be a semisimple Lie algebra. Then 

a) any Lie bialgebra $(\g,\de)$ is a coboundary one;

b) the quantization, $U_h(\g,r)$, of any coboundary Lie bialgebra $(\g,r)$ exists and
is isomorphic to $\Ug[[h]]$ as a topological $\C[[h]]$-algebra;

c) the comultiplication in $U_h(\g,r)$ has the form
\be{}
\label{comul}
\De_h(x)=F_h\De(x)F^{-1}_h, \qquad x\in\Ug,
\ee{}
 where $F_h\in \Ug^{\ot 2}[[h]]$ and can be chosen in the form
\be{}
\label{F}
F_h=1\ot 1+\hlf{h}r +o(h).
\ee{} 
Moreover,
\be{}\label{counit}
(\eps\ot 1)F_h=(1\ot\eps)F_h=1\ot 1,
\ee{}
where $\eps: U_h(\g,r)\to \C[[h]]$ is the counit in $U_h(\g,r)$, which coincides with
the natural extension of the counit $\Ug\to\C$.
\end{propn}

It follows from the coassociativity $\De_h$ that $F_h$  
satisfies  the
equation
\be{}
\label{eqF}
(F_h\t 1)\cdot(\De\t id)(F_h)=(1\t F_h)\cdot(id\t\De)(F_h)\cdot\Phi_h
\ee{}
for some invariant element $\Phi_h\in\Ug^{\ot 3}[[h]]$.

It follows from (\ref{eqF}) that if $F_h$ has the form (\ref{F}), then
the coefficient by $h$ in $\Phi_h$ vanishes. Moreover, the coefficient by $h^2$
is up to a factor the element $\ff$ from (\ref{cYB}), i.e.,  
\be{}\label{fPhi}
\Phi_h=1\t 1\t 1+h^2\ff+o(h^2).
\ee{}
In addition, it follows  from (\ref{eqF}) that
$\Phi_h$ satisfies the pentagon identity
\be{}
\label{pent}
(id^{\t 2}\t \De)(\Phi_h)\cdot(\De\t id^{\t 2})(\Phi_h)=
(1\t\Phi_h)\cdot(id\t\De\t id)(\Phi_h)\cdot(\Phi_h\t 1).
\ee{}

One can suppose that  an invariant element $\Phi_h$ of the form \bref{fPhi} satisfying
\bref{pent} is given in advance. Such $\Phi_h$ can be constructed independently, \cite{Dr2}. 
Then, for any element $r$
obeying \bref{cYB} one can find $F_h$ of the form \bref{F} 
such that \bref{eqF} holds. So, we may take such $F_h$ to
define comultiplication \bref{comul} in the quantum group $U_h(\g,r)$.
In the following we fix $\Phi_h$ satisfying the additional relation, \cite{DS2}:
\be{}\label{addPhi}
\Phi_h\Phi_h^{s}=1, 
\ee{}
where
$s$ is the antipod in $U(\g)$ defined by
$s(x)=-x$ for $x\in \g$, and $\Phi_h^{s}=(s\t s\t s)(\Phi_h)$. 

The Hopf algebra $U=U_h(\g,r)$ is quasitriangular with the universal R-matrix
\be{}\label{Rm}
\R=F_{21}e^{\frac{h}{2}{\bf t}}F^{-1}=1\ot 1+h(\frac{{\bf t}}{2}-r)+\cdots,
\ee{}
where ${\bf t}\in\Sym^2\g$ is an invariant element.
$\R$ satisfies the relations, \cite{Dr2}
\be{}\label{RDel}
\Delta'(x)=\R\Delta(x)\R^{-1}
\ee{}
and
\be{}\label{Rrel}
\Delta_1\R=\R_{13}\R_{23}\\ \notag
\Delta_2\R=\R_{13}\R_{12}.
\ee{}

Here subscripts indicate the position of the corresponding tensor factors;
for example, if $a=a_1\ot a_2$, then $a_{13}=a_1\ot 1\ot a_3$ and $a_{21}=a_2\ot a_1$.

\section{$\Uh$-invariant quantization}

Let $G$ be a semisimple connected complex Lie group with the Lie algebra  $\g$.
Let $M$ be a $G$-manifold and $\A=\A(M)$ a sheaf of 
functions on $M$. It may be the sheaf of analytic, smooth, or  
algebraic functions, dependingly of the type of $M$.
Then $\Ug$ acts on sections of $\A$, and the multiplication
in $\A$ is $\Ug$-invariant, i.e.,
\be{*}
b m(x\otimes y)=m\De(b)(x\t y)\quad\quad b\in\Ug,\ x,y\in \A.
\ee{*}
By a deformation quantization of $\A$ we mean a sheaf of associative
algebras, $\A_h$, which is 
equal to $\A[[h]]$ 
as a $\C[[h]]$-module, with multiplication in
$\A_h$ of  the form
\be{*}
m_h=\sum_{k=0}^\infty h^k m_k,
\ee{*}
where
$m_0=m$ is the usual commutative multiplication in $\A$
and $m_k$, $k>0$, are bidifferential operators vanishing on constants.
The action $\Ug$ on $\A$ is naturally extended to the action of
the $\C[[h]]$-algebra $\Ug[[h]]=\Uh$ on the $\C[[h]]$-module $\A[[h]]=\A_h$. 

We will also consider two parameter quantizations on $M$.
A two parameter quantization  of $\A$ is
an algebra $\A_{t,h}$ equal to $\A[[t,h]]$
as a $\C[[t,h]]$-module and having a multiplication of the form
\be{}\label{twom}
m_{t,h}=m_0+t m_1^\prime+h m_1^{\prime\prime}+o(t,h).
\ee{}

We are going to study  quantizations of $\A$ which are invariant under
the $\Uh$-action, i.e., under the comultiplication $\De_h$. The invariance means that
\be{}\label{conqinv}
b m_h(x\otimes y)=m_h\De_h(b)(x\t y)\quad\quad b\in\Ug,\ x,y\in \A,
\ee{}
where $m_h$ means one or two parameter multiplication.
Of course, for a two parameter $\Uh$-invariant multiplication, $m_{t,h}$,
the multiplication $m_t=m_{t,0}$ is $\Ug$-invariant.

We call $m_h$ satisfying \bref{conqinv} a $\Uh$-invariant quantization of
the $G$-manifold $M$.

\begin{defn} A $\C[[h]]$-linear map $\mu_h:\A_h\ot\A_h\to\A_h$ is called
a $\Phi_h$-{\em associative} multiplication if
\be{*}
\mu_h(\Phi_1x\otimes \mu_h(\Phi_2 y\t \Phi_3 z)))=\mu_h(\mu_h(x\t y)\t z)
\quad \mbox{ for } x,y,z\in \A,  
\ee{*}
where $\Phi_h=\Phi_1\ot\Phi_2\ot\Phi_3$ (summation implicit).

We say that the $\Phi_h$-associative multiplication 
$\mu_h=\sum_{k=0}^\infty h^k \mu_k$ gives a $\Phi_h$-associative
quantization of $\A$ if $\mu_0=m_0$,  the usual multiplication in $\A$,
and $\mu_k$, $k>0$, are bidifferential operators vanishing on constants.
\end{defn}

\begin{propn}\label{propo2.2} There is a natural one-to-one correspondence between
$\Uh$-invariant and $\Ug$-invariant $\Phi_h$-associative quantizations
of $\A$. Namely, if $\mu_h$ is a $\Ug$-invariant $\Phi_h$-associative 
multiplication in $\A[[h]]$, then
\be{}\label{cormul}
m_h=\mu_hF_h^{-1}
\ee{} 
gives a $\Uh$-invariant associative
multiplication in $\A[[h]]$.    
\end{propn}
\begin{proof} This follows immediately from (\ref{comul}) and (\ref{eqF}).
This follows also from the categorical interpretation of
$\Phi_h$ and $F_h$, \cite{Dr2}, \cite{DGS}.
\end{proof}

This proposition shows that given a $\Ug$-invariant $\Phi_h$-associative 
quantization of $\A$, we can get the $\Uh$-invariant 
quantization of $\A$ for any quantum
group $\Uh$ from Proposition \ref{propo2.1} b) by applying
$F_h$ from (\ref{F}) to the $\Phi_h$-associative multiplication. 

\section{Poisson brackets associated with $\Uh$-invariant\\ quantization}
A skew-symmetric map $f:\A^{\ot 2}\to\A$ we call a
{\em bracket} if it satisfies the Leibniz rule:
$f(ab,c)=af(b,c)+f(a,c)b$ for $a,b,c\in\A$.
It is easy to see that any bracket is presented by a bivector
field on $M$. Further we will identify brackets and
bivector fields on $M$.

For an element $\psi\in\wedge^k\g$ we denote by 
$\psi_M$ the $k$-vector field on $M$ which is induced
by the action map $\g\to\Vect(M)$.

A bracket $f$ is a Poisson one if the Schouten bracket
$\[f,f\]$ is equal to zero.

\begin{defn}\label{defphi} A $G$-invariant bracket $f$ on $M$ we call a
$\ff$-{\em bracket} if 
\be{}\label{phibr}
\[f,f\]=-\ff_M,
\ee{} 
where $\ff\in\wedge^3\g$ is an invariant element.
\end{defn}

\begin{propn}\label{propo2.3}
Let $\A_h$ be a $\Ug$-invariant $\Phi_h$-associative
quantization with multiplication 
$\mu_h=m_0+h\mu_1+o(h)$, where $m_0$ is
the multiplication in $\A$.
Then the map $f:\A^{\ot 2}\to\A$, $f(a,b)=\mu_1(a,b)-\mu_1(b,a)$,
is a $\ff$-bracket for $\ff$ from (\ref{fPhi}).
\end{propn}

\begin{proof} A direct computation.
\end{proof}

\begin{cor}\label{coro2.3} Let $\A_h$ be a $\Uh$-invariant associative
quantization with multiplication $m_h=m_0+h m_1+o(h)$.
Then the corresponding Poisson bracket $p(a,b)=m_1(a,b)-m_1(b,a)$
has the form
\be{}\label{dopq}
p(a,b)=f(a,b)-r_M(a,b),
\ee{}
where $r$ is the r-matrix
corresponding to $\Uh$ and $f$ is a $\ff$-bracket
with $\ff=\[r,r\]$.
\end{cor}

\begin{proof} Follows from Proposition \ref{propo2.3} and \bref{cormul}.
\end{proof}

Now, let us consider a two parameter $\Uh$-invariant quantization \bref{twom}. 
With such a quantization one associates two Poisson brackets:
the bracket $v(a,b)=m^\prime_1(a,b)-m^\prime_1(b,a)$ along $t$,
and the bracket
$p(a,b)=m^{\prime\prime}_1(a,b)-m^{\prime\prime}_1(b,a)$ along $h$.
It is easy to check that $p$ and $v$ are compatible Poisson brackets,
i.e., their Schouten bracket $\[p,v\]$ is equal to zero. So, we have

\begin{cor}\label{coro2.4} Let $\A_{t,h}$ be a $\Uh$-invariant
associative quantization of the form \bref{twom}.
Then the Poisson bracket  $p(a,b)=m_1^{\prime\prime}(a,b)-m_1^{\prime\prime}(b,a)$
has the form
\be{*}
p(a,b)=f(a,b)-r_M(a,b),
\ee{*}
where $r$ is the r-matrix
corresponding to $\Uh$ and $f$ is a $\ff$-bracket
with $\ff=\[r,r\]$. The Poisson bracket 
$s(a,b)=m_1^\prime(a,b)-m_1^\prime(b,a)$ is $\Ug$-invariant
and compatible with $f$ (and thus with  $p$), i.e., $f$ satisfies (\ref{phibr}) and the additional
condition
\be{}\label{phibrcom}
\[f,s\]=0.
\ee{}
\end{cor}
\begin{proof} Similar to Corollary \bref{coro2.3}.
\end{proof}

So, studying the problem of $\Uh$-invariant quantizations on a $G$-manifold $M$
is devided into two parts: 1) finding $G$-invariant bivector fields, $f$, on $M$
satisfying \bref{phibr} and \bref{phibrcom} (if a $G$-invariant Poisson bracket
$s$ is given); 2) proving the existence of quantization of brackets \bref{dopq}. 

\section{Quantization on symmetric spaces and on $G$}

\subsection{Quantization on symmetric spaces}

Let $G$ be a semisimple connected Lie group. Let us recall that a symmentric space
over $G$ is a space of the form $M=G/H$, where $H$ is
a subgroup of $G$ such that $G^\sigma_0\subset H\subset G^\sigma$,
where $G^\sigma$ is the set of fixed points of $\sigma$, an involutive automorphism of $G$,  
and $G^\sigma_0$
is the identity component of $G^\sigma$. The automorphism $\sigma$ induces
an automorphism of the Lie algebra $\g$ of $G$ that
we also denote by $\sigma$.

Let $r$ be an r-matrix on $\g$ and $U_h(\g,r)$ the corresponding quantum group.
In \cite{DS1}, the following statement is proven.
\begin{thm}\label{thmsym}
Let $M$ be a symmetric space over a semisimple Lie group $G$, $\sigma$ the corresponding
involution of $\g$.
Let an r-matrix $r$ on $\g$ be such that the element $\[r,r\]=\ff$ is $\sigma$-invariant.
Then there exists a  $U_h(\g,r)$-invariant quantization on $M$.
\end{thm}

\begin{rem} Note that the $\sigma$-invariance of $\ff$ implies that $\ff_M=0$ 
and $r_M$ is a Poisson bracket. So, in the case of symmetric space $M$
the invariant part $f$ of the corresponding Poisson bracket \bref{dopq}
may be taken to be zero.
\end{rem}

\subsection{Quantization of two-sided action $G$ on itself}
Consider the group $G$ as a $G\times G$-manifold with
the two-sided action of $G$ on itself.  Then,
the group $G$ may be considered as a symmetric
space, $G=(G\times G)/H$, where $H$ is the diagonal.
The action of $G\times G$ on $G$ is: $(g_1,g_2)\tr g=g_1gg^{-1}_2$,
$(g_1,g_2)\in G\times G,\ g\in G$.
In this case $\sigma$ on $\g\oplus\g$, the Lie algebra of $G\times G$,
is the usual permutation: $\sigma(x,y)=(y,x)$.

Let $r_1$ and $r_2$ be two r-matrices on $\g$ such that
\be{}\label{eqrm}
\[r_1,r_1\]=\[r_2,r_2\]=\ff.
\ee{}
The quantum group  $U_h(\g,r_1)\ot U_h(\g,r_2)$ may be considered as
a quantization of the Lie bialgebra $\g\oplus\g$ with r-matrix $(r_1,r_2)$. 
Applying Theorem \ref{thmsym} we obtain

\begin{propn}\label{proptwoq}
Let $G$ be a semisimple Lie group.
Let $r_1$ and $r_2$ be r-matrices on its Lie algebra satisfying \bref{eqrm}.
Then there exists a $U_h(\g,r_1)\ot U_h(\g,r_2)$-invariant quantization of the two-sided
action $G$ on $G$.
\end{propn}

\begin{rem} It is clear that the presence of a $U_h(\g,r_1)\ot U_h(\g,r_2)$-action is
equivalent to the presence of two commuting actions: $\rho_1$ of $U_h(\g,r_1)\ot 1$ and
$\rho_2$ of $1\ot U_h(\g,r_2)$. So, the previous proposition states the existence
of the invariant quantization on $G$ with respect to $\rho_1$ and $\rho_2$ simultaneously.
\end{rem}

\begin{rem} The element $\Phi_h$ from \bref{pent} associated to
$U_h(\g,r_1)$ and $U_h(\g,r_2)$ induces in the obvious way the corresponding element
$\Phi_h\ot\Phi_h$
associated to $U_h(\g,r_1)\ot U_h(\g,r_2)$. It follows from \bref{addPhi}
that the image of $\Phi_h\ot\Phi_h$ by the $G\times G$-action on $G$
is equal to unity. This implies that the usual multiplication
in the function algebra on $G$ is $\Phi_h\ot\Phi_h$-associative.
So, accordingly to Proposition \bref{propo2.2}, the quantization from 
the conclusion of Proposition \ref{proptwoq}
may be given explicitly as 
\be{}\label{twoact}
m_h=m\rho_1(F_{h1}^{-1})\rho_2(F_{h2}^{-1}),
\ee{}
where $F_{h1}$ and $F_{h2}$ are elements from \bref{comul}
associated to $U_h(\g,r_1)$ and $U_h(\g,r_2)$, respectively.
\end{rem}
\begin{rem} It is easy to calculate that, up to a factor, the Poisson bracket on $G$
corresponding to the quantization from Proposition \ref{proptwoq}
has the form
\be{}\label{ftwo}
r_1^L+r_2^R,
\ee{}
where $r_i^L$ ( $r_i^R$ ) is the extension of $r_i$ as a left- (right-) invariant
bivector field on $G$.
If we take $r_1=r$, where $r$ is the Sklyanin-Drinfeld r-matrix \bref{SDr}
and $r_2=-r$, then \bref{ftwo} becomes
\be{}\label{fSD}
r^L-r^R,
\ee{}
which is the Sklyanin-Drinfeld Poisson bracket on $G$.
\end{rem}

\subsection{Quantization of $Ad$-action $G$ on itself}

Consider the group $G$ as a $G$-manifold with respect to
the $Ad$-action: $g_1\tr g=g_1gg_1^{-1}$. The corresponding
action $\Ug$ on $\A(G)$ is obtained by the embedding
$\Delta:\Ug\to\Ug\ot\Ug$ and the two-sided action 
$\Ug\ot\Ug$ on $\A(G)$ considered above. In the quantum case, given a
$\Uh\ot\Uh$-invariant quantization, $m_h$, of $\A(G)$,
the corresponding embedding
$\Delta_h:\Uh\to\Uh\ot\Uh$ 
does not give a $\Uh$-invariant quantization of the function algebra,
because this embedding is not a morphism of coalgebras.
Nevertheless, the multiplication $m_h$ can be modified to give a quantization
of the $Ad$-action.

Namely, the following lemma holds.
\begin{lemma}\label{lemcormul}
Let $\A$ be an algebra with $\Uh\ot\Uh$-invariant multiplication $m_h$.
Let the action $\Uh\ot\Uh$ on $\A$ is given by two commutative actions $\rho_1$ and
$\rho_2$ of $\Uh$. Let the action $\Uh$ on $\A$ is induced by the embedding
$\Delta_h:\Uh\to\Uh\ot\Uh$.
Then the multiplication
\be{}
\mu'_h(a\ot b)=\mu_h(\rho_2(\R_1)a\ot\rho_1(\R_2)b)
\ee{}
is associative and $\Uh$-invariant.
Here $\R=\R_1\ot\R_2$ is R-matrix of the quantum group $\Uh$.
\end{lemma}
\begin{proof} A direct computation using relations \bref{Rrel}.
\end{proof}

To construct $Ad$-invariant quantizations from two-sided $\Uh$-module
algebras one may  use the following

\begin{lemma}\label{lemcorcop}
Let $\A$ be an algebra with $\Uh^{cop}$-invariant multiplication $m_h$,
where $\Uh^{cop}$ is $\Uh$ with the opposite comultiplication.
Then the multiplication
\be{}
\mu'_h(a\ot b)=\mu(\R_1a\ot\R_2b)
\ee{}
is $\Uh$-invariant.
\end{lemma}
\begin{proof} The same as of Lemma \ref{lemcormul}.
\end{proof}

\begin{cor} There exists a $\Uh$-invariant quantization of the $Ad$-action $G$ on $G$.
\end{cor}

\begin{rem} 
Using \bref{twoact} and Lemma \ref{lemcormul}, it is easy to calculate
that, up to a factor, the invariant part of the Poisson bracket corresponding to the
quantization from the previous corollary is:
\be{}\label{Adbr}
f(a,b)=(t_1^L a)(t_2^R b)-(t_1^L b)(t_2^R a),
\ee{}
where  ${\bf t}=t_1\ot t_2$ is the invariant element of $Sym^2\g$.
Here $x^L$ ( $x^R$ ) is the extension of $x\in\g$ as a left- (right-) invariant vector
field on $G$. 
\end{rem}

\section{Quantization on $\g^*$}

Let $\g$ be a complex Lie algebra. 
Then, the symmetric algebra $S\g$ can be considered as a function (polynomial) algebra
on $\g^*$. The algebra $U(\g)$ is included in the family of algebras
$(S\g)_t=T(\g)[t]/J_t$, where $J_t$ is  the ideal generated by the elements
of the form $x\ot y-\sigma(x\ot y)-t[x,y]$, $x,y\in\g$, $\sigma$ is
the permutation. By the PBW theorem, $(S\g)_t$ is a free module
over $\C[t]$. 
We have $(S\g)_0=S\g$, so this
family of quadratic-linear algebras gives a $U(\g)$ invariant 
quantization of $S\g$ by the Lie bracket $s$.

It turns out that for $\g=sl(n)$ this picture can be extended to 
the quantum case, \cite{Do2}. Namely, 

\begin{thm}\label{thm6.1}
Let $\g=sl(n)$. Then, there exist deformations, $\sigma_h$ and $[\cdot,\cdot]_h$, of 
both the mappings
$\sigma$ and $[\cdot,\cdot]$ such that
the two parameter family of algebras
$(S\g)_{t,h}=T(\g)[[h]][t]/J_{t,h}$, where $J_{t,h}$ is the ideal generated 
by elements
of the form $x\ot y-\sigma_h(x\ot y)-t[x,y]_h$, $x,y\in\g$,
gives a $\Uh$-invariant quantization of the Lie bracket $s$ on $\g^*$.
\end{thm}

Using Corollary \ref{coro2.4}, let us describe the Poisson brackets corresponding to $(S\g)_{t,s}$.

\begin{propn}\label{prop6.3}
The pair of brackets corresponding to the quantization $(S\g)_{t,h}$ consists of
two compatible Poisson brackets: 

$s$ (along $t$) is the Lie bracket;

$p$ (along $h$) is a quadratic Poisson bracket of the form
$p=f-\{\cdot,\cdot\}_r$, where $f$ is the invariant quadratic bracket
which is defined by a unique up to a factor
invariant map 
\be{}\label{quadbr}
f:\wedge^2\g\to S^2\g,
\ee{}
and $\{\cdot,\cdot\}_r$
is the $r$-matrix bracket. Moreover, $\[s,f\]=0$ and
$ \[f,f\]=-\overline{\ff}$, where $\overline{\ff}$
has the form $\overline{\ff}(a,b,c)=[\ff_1,a][\ff_2,b][\ff_3,c]$,
and $\ff=\ff_1\wedge\ff_2\wedge\ff_3=\[r,r\]$. 
Recall that $\ff$ is a unique up to a factor invariant element of $\wedge^3\g$.
\end{propn}

Note that the bracket $f$ from the above proposition can be 
restricted to any orbit of $SL(n)$ in $sl(n)^*$. This follows from
the following general statement, \cite{Do1}.

\begin{propn}\label{prop6.4}
Let $G$ be a semisimple Lie group with its Lie algebra $\g$, 
$s=[\cdot,\cdot]$  the Poisson-Lie bracket on $\g^*$.
Let $f=\{\cdot,\cdot\}$ be an invariant bracket on an open set of $\g^*$
such that the Schouten
bracket $\[s,f\]$ is a three-vector field, $\psi$,  which can be restricted to an orbit $\O$
of $G$. 
Then $f$ can be restricted to $\O$.
\end{propn}
Note that an invariant bivector field defined on a part of an orbit is in fact
defined on the whole orbit.
 
\begin{rem} Using quantized Verma modules, one can prove that the family
$(S\g)_{t,h}$ can be restricted to any semisimple orbit to give a $\Uh$-invariant
quantization of the Kirillov bracket on it.
\end{rem}

\begin{rem} One can prove,  \cite{Do1}, that for a simple $\g\neq sl(n)$
there are no $\ff$-brackets on $\g^*$ compatible with the Poisson-Lie bracket.
Hence, in this case a two parameter $\Uh$-invariant quantization on $\g^*$ 
does not exist.
\end{rem}

It seems that for  a simple $\g\neq sl(n)$ there exist no $\ff$-brackets
on $\g^*$, which implies that on $\g^*$ even a one parameter $\Uh$-invariant
quantization does not exist. This statement could be easily derived from Conjecture \ref{conj}
at the end of the paper.

\begin{rem}
Note that on a neighborhood of zero in $\g$ a  $\ff$-bracket exists.
Indeed, due to the map $\exp:\g\to G$ the bracket \bref{Adbr} may be carried over
from $G$ to a neighborhood of zero in $\g$.

Note also that the bracket of the form \bref{Adbr} is correctly defined on $gl(n)$
(considered as an associative algebra) and being
restricted to $sl(n)$ gives a quadratic bracket which coincides with 
\bref{quadbr}.
\end{rem} 

\begin{rem}\label{remorb}
It follows from the previous remark that for any simple $\g$, 
a $\ff$-bracket exists on any
orbit in $\g^*$ (not necessarily semisimple).
Indeed, multiplying a given orbit by a constant one can suppose that
it goes through the neighborhood of zero where the bracket \bref{Adbr}
is defined. Now, the statement follows from Proposition \ref{prop6.4}.

Since bracket \bref{quadbr} is compatible
with the Poisson-Lie bracket, it follows that for $\g=sl(n)$ on any orbit
in $\g^*$ there exists a $\ff$-bracket compatible with the Kirillov bracket.

As we will see in Subsection \ref{ss7.2}, in case $\g\neq sl(n)$ not all
orbits (even semisimple ones) admit 
$\ff$-brackets compatible with the Kirillov ones.
\end{rem} 

\section{Quantization on semisimple orbits}

Let $G$ be a complex connected simple Lie group with the Lie
algebra $\g$.
Let $\lf$ be a Levi subalgebra of $\g$, the Levi factor of a parabolic subalgebra.
Let $L$ be a Lie subgroup of $G$ with Lie algebra $\lf$. Such a subgroup is called
a Levi subgroup.
It is known that $L$ is a closed connected subgroup. Denote $M=G/L$ and
let $o\in M$ be the image of the unity by the natural projection $G\to M$.
Then $L$ is the stabilizer of $o$. It  is known, that $M$ may be
realized as a semisimple orbit of $G$ in the coadjoint representation $\g^*$. 
Conversely, any semisimple orbit in $\g^*$ is a quotient of $G$ by a Levi subgroup.  
We call {\em the rank} of $M$ the dimension of the center of $\lf$. So, if
$M$ is a maximal orbit, i.e., $\lf$ is equal to the Cartan subalgebra, then the rank of $M$
is equal to the rank of $\g$.

\subsection{One parameter quantization}

In \cite{Do3} the following statement is proven.
\begin{thm} Let $G$ be a simple Lie group, $\g$ its Lie algebra, $M$ a semisimple
orbit in $\g^*$. Then 

a) The set of all $\ff$-brackets on $M$ form an affine nonsingular algebraic manifold, $X_M$,
of dimension $rank(M)$.

b) There exists an analytic universal $\Uh$-invariant family of multiplications 
$\mu_{f,h}$ on $\A(M)$
of the form
\be{}\label{aspr}
\mu_{f,h}(a,b)=ab +(h/2) (f(a,b)-r_M(a,b)) +\sum_{n\geq 2}h^n \mu_{f,n}(a,b), \quad f\in X.
\ee{}
The universality means that
for any  $\Uh$-invariant multiplication,
$m_h$, there exists a formal path in $X_M$, $\psi(h)$, such that
$m_h$ is equivalent to $\mu_{\psi(h),h}$, and
multiplications corresponding to different paths are
not equivalent.
\end{thm}

\begin{rem}
Let us consider $M$ as an abstract manifold. Then all nondegenerate $G$-invariant
Poisson brackets on it appear as restrictions of the Poisson-Lie bracket
by embeddings $M$ into $\g^*$ as orbits. In the quantum case, when $\g=sl(n)$,
a dense open set of $\ff$-brackets in $X_M$ can be obtained in the same way,
as restrictions of the quadratic bracket \bref{quadbr}.
But in general there are $\ff$-brackets on $M$ which can not be obtained in this way.
For example, the Sklyanin-Drinfeld bracket \bref{fSD} can be carried over from $G$
to $M$ considered as a quotient of $G$, but the obtained bracket can not be
realized by embedding $M$ into $\g^*$ as an orbit. 

Note that also for $\g\neq sl(n)$ the essential part of $\ff$-brackets on $M$
may be obtained with the help of embeddings $M$ into $\g^*$ 
as restrictions of the bracket \bref{Adbr} carried over to an open set of $\g^*$ 
(see Remark \ref{remorb}).
\end{rem}

\subsection{Two parameter quantization}
\label{ss7.2}
It turns out that $\ff$-brackets compatible with the Kirillov bracket exist
not on any orbits in $\g^*$.

\begin{defn} 
Let $M$ be an orbit in $\g^*$ (not necessarily semisimple).
We call $M$ a {\em good} orbit, if there exists a $\ff$-bracket on it
compatible with the Kirillov bracket (we call such a bracket a {\em good}
bracket).
\end{defn}

The following proposition, \cite{DGS}, gives a classification of 
good semisimple orbits.

\begin{propn}\label{prop7.7}

a) For $\g$ of type $A_n$ all semisimple orbits are good.

b) For  all other $\g$, a semisimple orbit $M$ is good if and only if
$M=G/L$. Here $L$ is a Levi subgroup whose Lie algebra is generated
by the Cartan subalgebra and root vectors $X_{\pm\alpha}$, where
$\alpha$ runs all simple roots away of
one or two roots which appear   
in the representation of the maximal root with coefficient 1.

c) The good brackets $f$ on a good orbit, $M$, form a one-dimensional variety
isomorphic to $\C$:
all such brackets have the form 
\be{*}
\pm f_0+c\cdot s_M,
\ee{*}
where $s_M$ is the Kirillov bracket on $M$, $c\in \C$, and $f_0$ is a fixed good bracket.
\end{propn}
So, if $\g\neq sl(n)$, the good orbits are either Hermitian symmetric spaces
or, for $\g$ of types $D_n$ and $E_6$, bisymmetric spaces. 

It follows from Proposition \ref{prop7.7} c) that for a good orbit the
family 
\be{}\label{fambr}
h(f_0-r_M)+ts_M,
\ee{} 
is uniquely defined up to a linear change of parameters $(h,t)$.

\begin{propn} Let $M$ be a good semisimple orbit. Then the family \bref{fambr}
can be quantized, i.e. 
there exists a two parameter $\Uh$-invariant multiplication in $\A(M)$ of the form
\be{*}
m_{t,h}(a,b)=ab+\frac{1}{2}(h(f_0-r_M)+ts_M)+\sum_{k+l\geq 2}h^kt^lm_{k,l}(a,b).
\ee{*}
\end{propn}
The theorem is proven in \cite{Do3}.

\subsection{Problems}
An open question is to investigate the  problem of one and two parameter quantization
for non-semisimple orbits, in particular, for nilpotent ones. In this connection,
I would like to mention the following problems which seem to be
purely from the theory of (simple) Lie algebras:

1. Give a classification of good orbits in case $\g\neq sl(n)$ (see Remark \ref{remorb}).

2. Give a classification of $\ff$- and good brackets for non-semisimple orbits
(see the same Remark).

3. Prove the following 
\begin{conjecture}\label{conj}
Let $\g$ be a simple Lie algebra. Then all invariant Poisson
brackets on the polynomial algebra on $\g^*$, $S\g$, 
have the form $b[\cdot,\cdot]$, where
$[\cdot,\cdot]$ is the Poisson-Lie bracket and $b$ an invariant element of $S\g$.
\end{conjecture}

\end{document}